\documentstyle[11pt]{article}

\newtheorem{lemma}{Lemma}[section]
\newtheorem{theorem}{Theorem}[section]
\newtheorem{corollary}{Corollary}[section]

\def\blemma{\begin{lemma}\sl{}\def\elemma{\end{lemma}}}
\def\btheorem{\begin{theorem}\sl{}\def\etheorem{\end{theorem}}}

\def\beqlb{\begin{eqnarray}}\def\eeqlb{\end{eqnarray}}
\def\beqnn{\begin{eqnarray*}}\def\eeqnn{\end{eqnarray*}}

\def\qed{\quad$\Box$\medskip}

\def\lan{\langle}\def\ran{\rangle}

\def\E{\mbox{\boldmath$E$}}

\def\IR{I\!\!R}

\begin{document}

\

\noindent{Published in: {\it Potential Analysis} {\bf 20} (2004),
285--302.}

\bigskip\bigskip

\noindent{\Large\bf Non-differentiable Skew Convolution
Semigroups}

\medskip
\noindent{\Large\bf and Related Ornstein-Uhlenbeck Processes}

\bigskip
\noindent{DONALD A. DAWSON\,\footnote{Supported by an NSERC
Research Grant and a Max Planck Award.}}

\smallskip
\noindent{\it School of Mathematics and Statistics, Carleton
University, 1125 Colonel By Drive,}

\noindent{\it Ottawa, Canada K1S 5B6 (e-mail:
ddawson@math.carleton.ca)}

\bigskip
\noindent{ZENGHU LI\,\footnote{Supported by the NSFC
(No.\,10131040 and No.\,10121101).}}

\smallskip
\noindent{\it Department of Mathematics, Beijing Normal
University, Beijing 100875, P.R. China}

\noindent{\it (e-mail: lizh@email.bnu.edu.cn)}

\bigskip\bigskip

\noindent{\bf Abstract:} It is proved that a general
non-differentiable skew convolution semigroup associated with a
strongly continuous semigroup of linear operators on a real
separable Hilbert space can be extended to a differentiable one on
the entrance space of the linear semigroup. A c\`adl\`ag strong
Markov process on an enlargement of the entrance space is
constructed from which we obtain a realization of the
corresponding Ornstein-Uhlenbeck process. Some explicit
characterizations of the entrance spaces for special linear
semigroups are given.

\bigskip


\bigskip


\bigskip\bigskip

\section{Introduction}
\setcounter{equation}{0}

Suppose that $(S,+)$ is a Hausdorff topological semigroup and
$(Q_t)_{t\ge0}$ is a transition semigroup on $S$ satisfying
 \beqlb\label{1.1}
Q_t(x_1+x_2,\cdot) = Q_t(x_1,\cdot) * Q_t(x_2,\cdot),
\quad t\ge0, x_1,x_2\in S,
 \eeqlb
where ``$*$'' denotes the convolution operation. A family of
probability measures $(\mu_t)_{t\ge0}$ on $S$ is called a {\it
skew convolution semigroup} (SC-semigroup) associated with
$(Q_t)_{t\ge0}$ if it satisfies
 \beqlb\label{1.2}
\mu_{r+t} = (\mu_rQ_t) * \mu_t,
\quad r,t\ge0.
 \eeqlb
This equation is of interest since it is satisfied if and only if
 \beqlb\label{1.3}
Q_t^\mu (x,\cdot) := Q_t(x,\cdot) * \mu_t(\cdot),
\quad t\ge0, x\in S,
 \eeqlb
defines another transition semigroup $(Q_t^\mu)_{t\ge0}$ on $S$.
(Note that (\ref{1.1}) implies $(\mu*\nu)Q_t = (\mu Q_t) * (\nu
Q_t)$ for probability measures $\mu$ and $\nu$ on $S$.) This fact
was first observed in \cite{L95/6, L96} when $S= M(E)$
is the space of all finite Borel measures on a metrizable space
$E$; see also \cite[Theorem~2.1]{L02}. In that case,
$(Q_t)_{t\ge0}$ corresponds to a measure-valued branching process
and $(Q_t^\mu)_{t\ge0}$ corresponds to an immigration process.

In this work, we shall consider the formulation in another
special situation, where $S=H$ is a real separable Hilbert space
and $Q_t(x,\cdot)\equiv \delta_{T_tx}$ for a strongly continuous
semigroup of bounded linear operators $(T_t)_{t\ge0}$ on $H$. In
this case, we can rewrite (\ref{1.2}) as
 \beqlb\label{1.4}
\mu_{r+t} = (T_t\mu_r) * \mu_t,
\quad r,t\ge0,
 \eeqlb
and the transition semigroup $(Q^\mu_t)_{t\ge0}$ is given by
 \beqlb\label{1.5}
Q^\mu_tf(x) := \int_H f(T_tx+y)\mu_t(dy),
\quad x\in H,f\in B(H),
 \eeqlb
where $B(H)$ denotes the totality of bounded Borel measurable
functions on $H$. The semigroup $(Q^\mu_t)_{t\ge0}$ defined by
(\ref{1.5}) is called a {\it generalized Mehler semigroup}
associated with $(T_t)_{t\ge0}$, which corresponds to a
generalized Ornstein-Uhlenbeck process (OU-process). This
definition of the generalized Mehler semigroup was given by
Bogachev {\it et al} \cite{BRS96}. They also gave a
characterization for the SC-semigroup $(\mu_t)_{t\ge0}$ under the
assumption that the function $t \mapsto \hat\mu_t(a)$ is
differentiable at $t=0$, where $\hat \mu_t(a)$ denotes the
characteristic functional of $\mu_t$. It is known that for a
general SC-semigroup $(\mu_t)_{t\ge0}$ defined by (\ref{1.4}) the
function $t \mapsto \hat\mu_t(a)$ is not necessarily
differentiable at $t=0$; see e.g.\ \cite{DLSS03, N00, SS01}. A
simple and nice necessary and sufficient condition for an
SC-semigroup to be differentiable was given in \cite{N00} in the
setting of cylindrical probability measures. In \cite{BRS96} it
was shown that a differentiable cylindrical Gaussian SC-semigroup
can be extended into a real Gaussian SC-semigroup in an
enlargement of $H$ and the corresponding OU-process was
constructed as the strong solution to a stochastic differential
equation. Those results were extended to the general non-Gaussian
case in \cite{FR00}. A characterization for general SC-semigroups
$(\mu_t)_{t\ge0}$ was given in \cite{DLSS03}, where it was also
observed that the OU-processes corresponding to a
non-differentiable SC-semigroup usually have no right continuous
realizations. This property is similar to that of the immigration
processes studied in \cite{L96, L98, L02} and represents a
departure from the theory of well-studied classes of OU-processes
in \cite{BRS96, FR00}.

The main purpose of this paper is to study the construction of
OU-processes corresponding to non-differentiable SC-semigroups. We
shall see that, under a moment assumption, a general SC-semigroup
can be decomposed as the convolution of a centered SC-semigroup
and a degenerate one. For this reason, we shall only consider
centered SC-semigroups. In section 2, we derive from the results
of Dawson {\it et al} \cite{DLSS03} that each centered
SC-semigroup is uniquely determined by an infinitely divisible
probability measure on the entrance space $\tilde H$ for the
semigroup $(T_t)_{t\ge0}$, which is an enlargement of $H$. In
section 3, it is shown that a general non-differentiable centered
SC-semigroup can always be extended to a differentiable one on the
entrance space $\tilde H$. In section 4, we use a modification of
the argument of Fuhrman and R\"ockner \cite{FR00} to construct a
c\`adl\`ag and strong Markov OU-process $\{\bar X_t: t\ge0\}$ on a
further extension $\bar H$ of $\tilde H$. We also show that, if
$\bar X_0\in H$, then $\bar X_t\in H$ almost surely for every
$t\ge0$ and $\{1_H(\bar X_t)\bar X_t: t\ge0\}$ is an OU-process
with transition semigroup $(Q^\mu_t)_{t\ge0}$ defined by
(\ref{1.5}). Those results provide an approach to the study of
non-differentiable generalized Mehler semigroups with which one
can reduce some of their analysis to the existing framework of
\cite{BRS96} and \cite{FR00}. However, this approach should not
convince the reader that non-differentiable generalized Mehler
semigroups do not bear particular consideration on their own. In
fact, there are some cases where the natural state space of the
OU-process is $H$ and the introduction of $\tilde H$ and $\bar H$
seems unnatural and artificial. For example, an OU-process on
$L^2(0,\infty)$ with non-differentiable SC-semigroup represents
the fluctuation density of a catalytic branching processes with
immigration; see \cite{DLSS03}. In this case, it is rather
unnatural to take $\bar L^2(0,\infty)$ as the state space. We
provide some explicit characterization for the non-negative
elements of $\tilde L^2(\IR^d)$ and $\tilde L^2(0,\infty)$ in
section 5. The explicit characterization for all elements of $\bar
L^2(\IR^d)$ and $\bar L^2(0,\infty)$ seems much more
sophisticated.

\section{Non-differentiable semigroups}
\setcounter{equation}{0}

Suppose that $H$ is a real separable Hilbert space with dual space
$H^*$ and $(T_t)_{t\ge0}$ is a strongly continuous semigroup of
linear operators on $H$ with dual $(T_t^*)_{t\ge0}$. Let
$(\mu_t)_{t\ge0}$ be an SC-semigroup defined by (\ref{1.4})
satisfying the moment condition
 \beqlb\label{2.1}
\int_{H^\circ} \|x\|^2 \mu_t(dx) < \infty,
\quad t\ge 0,
 \eeqlb
where $H^\circ = H\setminus\{0\}$. Then we may define an
$H$-valued path $(b_t)_{t\ge0}$ by Bochner integrals
 \beqnn
b_t := \int_{H^\circ} x \mu_t(dx),
\quad t\ge 0,
 \eeqnn
and define $\mu^c_t = \delta_{-b_t}*\mu_t$. It is easy to check
that both $(\delta_{b_t})_{t\ge0}$ and $(\mu^c_t)_{t\ge0}$ are
SC-semigroups associated with $(T_t)_{t\ge0}$ and $\mu_t =
\mu^c_t*\delta_{b_t}$. That is, under the moment assumption, a
general SC-semigroup can be decomposed as the convolution of a
centered SC-semigroup and a degenerate one. For this reason, we
shall only discuss centered SC-semigroups in the sequel.

Since $(T_t)_{t\ge0}$ is strongly continuous, there are constants
$c_0\ge 0$ and $b_0\ge0$ such that $\|T_t\| \le c_0 e^{b_0t}$. Let
$(U_\alpha)_{\alpha>b_0}$ denote the resolvent of $(T_t)_{t\ge 0}$
and let $A$ denote its generator with domain $D(A) = U_\alpha H
\subset H$. An $H$-valued path $\tilde x = \{\tilde x(s): s>0\}$
is called an {\it entrance path} for the semigroup $(T_t)_{t\ge0}$
if it satisfies $\tilde x(r+t) = T_t\tilde x(r)$ for all $r,t>0$.
Let $E$ denote the set of all entrance paths for $(T_t)_{t\ge0}$.
We say $\tilde x\in E$ is {\it closable} if there is an element
$\tilde x(0)\in H$ such that $\tilde x(s) = T_s\tilde x(0)$ for
all $s>0$; and we say it is {\it locally square integrable} if
 \beqlb\label{2.2}
\int_0^l \|\tilde x(s)\|^2 ds < \infty
 \eeqlb
for some $l>0$.

\blemma\label{l2.1}
For any $\tilde x\in E$, (\ref{2.2}) holds for some $l>0$ if and
only if it holds for all $l>0$; and if and only if
 \beqlb\label{2.3}
\int_0^\infty e^{-2bs}\|\tilde x(s)\|^2 ds < \infty
 \eeqlb
for all $b>b_0$.
\elemma

{\it Proof.} Suppose that (\ref{2.2}) holds for some $l_0>0$. Let
$l>0$ and let $n\ge1$ be an integer such that $nl_0\ge l$. Then
 \beqnn
\int_0^l \|\tilde x(s)\|^2 ds
&\le&
\int_0^{nl_0} \|\tilde x(s)\|^2 ds  \\
&=&
\sum_{k=0}^{n-1}\int_0^{l_0} \|T_{kl_0}\tilde x(s)\|^2 ds  \\
&\le&
\sum_{k=0}^{n-1}c_0^2e^{2kb_0l_0}\int_0^{l_0}
\|\tilde x(s)\|^2 ds  \\
&<&\infty.
 \eeqnn
Thus (\ref{2.2}) holds for all $l>0$. On the other hand, for
any $b>b_0$,
 \beqnn
\int_0^\infty e^{-2bs}\|\tilde x(s)\|^2 ds
&=&
\sum_{k=0}^\infty e^{-2kbl_0}\int_0^{l_0} e^{-2bs}\|T_{kl_0}
\tilde x(s)\|^2 ds  \\
&\le&
\sum_{k=0}^\infty c_0^2e^{-2k(b-b_0)l_0}
\int_0^{l_0} e^{-2bs}\|\tilde x(s)\|^2 ds  \\
&<&\infty.
 \eeqnn
That is, (\ref{2.3}) holds for all $b>b_0$. The remaining assertions
are obvious.
\qed

Let $\tilde H$ denote the set of all locally square integrable
entrance paths for $(T_t)_{t\ge0}$. We shall call $\tilde H$ the
{\it entrance space} for $(T_t)_{t\ge0}$. For any fixed $b>b_0$,
we may define an inner product on $\tilde H$ by
 \beqlb\label{2.4}
\lan\tilde x,\tilde y\ran_\sim
:=
\int_0^\infty e^{-2bs}\lan \tilde x(s),\tilde y(s)\ran ds,
\quad \tilde x, \tilde y\in \tilde H.
 \eeqlb
Let $\|\cdot\|_\sim$ denote the norm induced by this inner product.
The proof of the following result was suggested to us by W. Sun.

\blemma\label{l2.2}
The normed space $(\tilde H,\|\cdot\|_\sim)$ is complete, so
$(\tilde H,\lan\cdot,\cdot\ran_\sim)$ is a Hilbert space.
\elemma

{\it Proof.} Suppose $\{\tilde x_n\} \subset \tilde H$ is a Cauchy
sequence under the norm $\|\cdot\|_\sim$, that is,
 \beqnn
\|\tilde x_n - \tilde x_m\|_\sim
=
\int_0^\infty e^{-2bs}\|\tilde x_n(s) - \tilde x_m(s)\|^2ds
\to 0
 \eeqnn
as $m,n \to \infty$. For each $t>0$,
 \beqnn
\|\tilde x_n(t) - \tilde x_m(t)\|^2
&=&
t^{-1}\int_0^t \|\tilde x_n(t) - \tilde x_m(t)\|^2 ds \\
&=&
t^{-1}\int_0^t \|T_{t-s}(\tilde x_n(s)
- \tilde x_m(s))\|^2 ds \\
&\le&
c_0^2t^{-1}e^{2bt}\int_0^t e^{-2bs}\|\tilde x_n(s)
- \tilde x_m(s)\|^2 ds.
 \eeqnn
Then the limit $\tilde x(t) = \lim_{n\to \infty} \tilde x_n(t)$
exists in $H$. Since $T_s$ is a continuous operator on $H$, for
$s>0$,
 \beqnn
T_s\tilde x(t)
= \lim_{n\to \infty} T_s\tilde x_n(t)
= \lim_{n\to \infty} \tilde x_n(t+s)
= \tilde x(t+s),
 \eeqnn
that is, $\tilde x = \{\tilde x(t): t>0\}$ is an entrance path
for $(T_t)_{t\ge0}$. For $\varepsilon>0$, choose large enough
$N\ge1$ such that
 \beqnn
\|\tilde x_n - \tilde x_m\|_\sim^2
=
\int_0^\infty e^{-2bs}\|\tilde x_n(s) - \tilde x_m(s)\|^2 ds
<\varepsilon
 \eeqnn
for $m,n\ge N$. By Fatou's lemma we get
 \beqnn
\int_0^\infty e^{-2bs}\|\tilde x_n(s) - \tilde x(s)\|^2ds
\le
\liminf_{m\to \infty} \int_0^\infty e^{-2bs} \|\tilde x_n(s)
- \tilde x_m(s)\|^2 ds
\le
\varepsilon.
 \eeqnn
It follows that
 \beqnn
\int_0^\infty e^{-2bs}\|\tilde x(s)\|^2ds
\le
\int_0^\infty e^{-2bs}\|\tilde x_n(s)\|^2ds
+ \int_0^\infty e^{-2bs} \|\tilde x_n(s) - \tilde x(s)\|^2 ds
< \infty.
 \eeqnn
Then $\tilde x\in \tilde H$ and $\lim_{n\to \infty} \| x_n
- x\|_\sim^2 =0$. \qed

\blemma\label{l2.3} The map $J: x \mapsto \{T_sx: s>0\}$ from
$(H,\|\cdot\|)$ to $(\tilde H,\|\cdot\|_\sim)$ is a continuous
dense embedding and hence $(\tilde H,\|\cdot\|_\sim)$ is
separable. \elemma

{\it Proof.} Since $x = \lim_{t\to 0^+} T_tx$, the map $J: x
\mapsto \{T_sx: s>0\}$ is injective. If $\lim_{n\to \infty} x_n =
x \in H$, then
 \beqnn
\int_0^\infty e^{-2bs}\|T_sx_n - T_sx\|^2ds
\le
c_0^2\|x_n - x\|^2\cdot \int_0^\infty e^{-2(b-b_0)s}ds
\to 0
 \eeqnn
as $n\to \infty$. Thus $J$ is a continuous embedding. For an
arbitrary $\tilde x \in \tilde H$ we have
 \beqnn
\| J\tilde x(t) - \tilde x\|_\sim^2
&=&
\int_0^\infty e^{-2bs}\|T_t\tilde x(s)
- \tilde x(s)\|^2 ds   \\
&=&
\int_0^r e^{-2bs}\|T_t\tilde x(s) - \tilde x(s)\|^2 ds
+ \int_r^\infty e^{-2bs}\|T_{s-r}[T_t\tilde x(r)
- \tilde x(r)]\|^2 ds   \\
&\le&
2(c_0^2e^{2b_0t} + 1)\int_0^r e^{-2bs} \|\tilde x(s)\|^2 ds \\
& & + c_0^2e^{-2b_0r}\|T_t\tilde x(r)
- \tilde x(r)\|^2 \int_r^\infty e^{-2(b-b_0)s}ds.
 \eeqnn
Observe that the first integral on the right hand side goes to
zero as $r\to 0^+$ and for fixed $r>0$ the second term goes to
zero as $t\to 0^+$. Then we have $\| J\tilde x(t) - \tilde
x\|_\sim \to 0$ as $t\to 0$, and $JH$ is dense in $\tilde H$.
Since $H$ is separable, so is $\tilde H$. \qed

\btheorem\label{t2.1} A family $(\mu_t)_{t\ge0}$ of centered
probability measures on $H$ satisfying (\ref{2.1}) is an
SC-semigroup associated with $(T_t)_{t\ge0}$ if and only if its
characteristic functionals are given by
 \beqlb\label{2.5}
\hat\mu_t(a)
=
\exp\bigg\{\int_0^t \log\hat\nu_s(a) ds\bigg\},
\quad t\ge0, a\in H^*,
 \eeqlb
where $(\nu_s)_{s>0}$ is a family of centered infinitely divisible
probability measures on $H$ satisfying $\nu_{r+t} = T_t\nu_r$ for
all $r,t>0$ and
 \beqlb\label{2.6}
\int_0^tds\int_{H^\circ} \|x\|^2 \nu_s(dx) < \infty,
\quad t\ge 0,
 \eeqlb
and $\log\hat\nu_s(\cdot)$ denotes the unique continuous function
on $H^*$ with $\log\hat\nu_s(0) =0$ and $\hat\nu_s(a) =
\exp\{\log\hat\nu_s(a)\}$ for all $a\in H^*$. \etheorem

{\it Proof.} It is well-known that the second moment of a centered
infinitely divisible probability measure only involves the
Gaussian covariance operator and the L\'evy measure. If the
centered probability measures $(\mu_t)_{t\ge0}$ and
$(\nu_s)_{s>0}$ are related by (\ref{2.5}), the Gaussian
covariance operators and L\'evy measures of $(\mu_t)_{t\ge0}$ can
be represented as integrals of those of $(\nu_s)_{s>0}$. This
observation yields that
 \beqnn
\int_{H^\circ} \lan x,a\ran ^2 \mu_t(dx)
=
\int_0^tds\int_{H^\circ} \lan x,a\ran ^2 \nu_s(dx),
\quad t\ge 0, a\in H^*.
 \eeqnn
Let $\{e_n: n=1,2,\dots\}$ be an orthonormal basis of $H=H^*$.
Applying the above equation to each $e_n$ and taking the summation
we see
 \beqlb\label{2.7}
\int_{H^\circ} \|x\|^2 \mu_t(dx)
=
\int_0^tds\int_{H^\circ} \|x\|^2 \nu_s(dx),
\quad t\ge 0.
 \eeqlb
Thus conditions (\ref{2.1}) and (\ref{2.6}) are equivalent for the
probability measures $(\mu_t)_{t\ge0}$ and $(\nu_s)_{s>0}$ related
by (\ref{2.5}). Suppose $(\mu_t)_{t\ge0}$ is given by (\ref{2.5})
with the centered infinitely divisible probabilities
$(\nu_s)_{s>0}$ satisfying $\nu_{r+t} = T_t\nu_r$ for all $r,t>0$.
Then $(\mu_t)_{t\ge0}$ is a centered SC-semigroup by
\cite[Theorem~2.3]{DLSS03}. Conversely, by \cite[Theorems~2.1 and
2.2]{DLSS03} any SC-semigroup $(\mu_t)_{t\ge0}$ has the expression
(\ref{2.5}) up to the convolution of a family of degenerate
probability measures $(\delta_{b_t})_{t\ge0}$. If
$(\mu_t)_{t\ge0}$ is a centered SC-semigroup, we must have $b_t=0$
for all $t\ge0$. \qed

\btheorem\label{t2.2} A family $(\mu_t)_{t\ge0}$ of centered
probability measures on $H$ satisfying (\ref{2.1}) is an
SC-semigroup associated with $(T_t)_{t\ge0}$ if and only if its
characteristic functionals are given by
 \beqlb\label{2.8}
\hat\mu_t(a)
=
\exp\bigg\{-\int_0^t\psi_s(a) ds\bigg\},
\quad t\ge0, a\in H^*,
 \eeqlb
where $\psi_s(\cdot)$ denotes the unique continuous function on
$H^*$ with $\psi_s(0) =0$ and
 \beqlb\label{2.9}
\exp\{-\psi_s(a)\}
=
\int_{\tilde H} e^{i\lan \tilde x(s),a\ran }\lambda_0 (d\tilde x),
\quad s>0, a\in H^*,
 \eeqlb
where $\lambda_0$ is a centered infinitely divisible probability
measure on $\tilde H$ satisfying
 \beqlb\label{2.10}
\int_{\tilde H} \|\tilde x\|_\sim^2 \lambda_0(d\tilde x)
< \infty.
 \eeqlb
\etheorem

{\it Proof.} Let $(\nu_s)_{s>0}$ be given as in Theorem
\ref{t2.1}. In the terminology of Markov processes,
$(\nu_s)_{s>0}$ is a probability entrance law for the Markov
process $\{T_tx: t\ge0\}$ with deterministic motion. Let $E_0 =
H^{(0,\infty)}$ be the totality of paths $\{w(t): t>0\}$ from
$(0,\infty)$ to $H$. We endow $E_0$ with the $\sigma$-algebra
${\cal E}_0$ generated by the maps $w \mapsto w(s)$, $s>0$. By
Kolmogorov's existence theorem, there is a unique probability
measure $\lambda_0$ on $E_0$ so that $\{w(t): t>0\}$ under
$\lambda_0$ is a Markov process with the same transition semigroup
as the process $\{T_tx: t\ge0\}$ and $\nu_s$ is the image of
$\lambda_0$ under $w \mapsto w(s)$; see e.g.\ Sharpe \cite[p.6]{S88}.
Because of the special deterministic motion mechanism of $\{T_tx:
t\ge0\}$ we may assume that $\lambda_0$ is supported by the
entrance paths $E$. Let ${\cal E}_0(E)$ and ${\cal E}_0(\tilde H)$
denote respectively the traces of ${\cal E}_0$ on $E$ and $\tilde
H$. Since $w \mapsto \|w(s)\|^2$ is clearly a non-negative ${\cal
E}_0(E)$-measurable function on $E$,
 \beqnn
w \mapsto \|w\|_\sim := \int_0^\infty e^{-2bs}\|w(s)\|^2ds
 \eeqnn
is an ${\cal E}_0(E)$-measurable function on $E$ taking values
in $[0,\infty]$. It is also easy to check that ${\cal E}_0(\tilde
H)$ coincides with the Borel $\sigma$-algebra ${\cal B}(\tilde H)$
induced by the norm $\|\cdot\|_\sim$. Since $(\nu_s)_{s>0}$
satisfies (\ref{2.6}), we have
 \beqnn
\int_E \|w\|_\sim^2 \lambda_0(dw)
&=&
\int_E \lambda_0(dw) \int_0^\infty e^{-2bs}\|w(s)\|^2 ds  \\
&=&
\int_0^\infty ds \int_H e^{-2bs}\|x\|^2 \nu_s(dx)    \\
&=&
\sum_{n=0}^\infty\int_0^1 ds \int_H e^{-2b(n+s)}
\|T_nx\|^2 \nu_s(dx)   \\
&\le&
c_0^2\sum_{n=0}^\infty e^{-2(b-b_0)n}\int_0^1 ds
\int_H e^{-2bs}\|x\|^2 \nu_s(dx)   \\
&<& \infty,
 \eeqnn
so $\lambda_0$ is supported by $\tilde H$ and (\ref{2.10})
holds. The infinite divisibility of $\lambda_0$ follows
immediately from that of $\nu_s$. \qed

\section{Differentiable extensions}
\setcounter{equation}{0}

For a general SC-semigroup given by Theorem \ref{t2.2}, the
function $t\mapsto \hat\mu_t(a)$ is not necessarily differentiable
at $t=0$. However, if $\nu_0$ is a centered infinitely divisible
probability measure on $H$ satisfying
 \beqlb\label{3.1}
\int_{H^\circ} \|x\|^2 \nu_0(dx) < \infty,
 \eeqlb
then
 \beqlb\label{3.2}
\hat\mu_t(a)
=
\exp\bigg\{\int_0^t \log\hat\nu_0(T^*_sa) ds\bigg\},
\quad t\ge0, a\in H^*,
 \eeqlb
defines a centered SC-semigroup $(\mu_t)_{t\ge0}$ such that
$t\mapsto \hat\mu_t(a)$ is differentiable at $t=0$ for all $a\in
H^*$. In the sequel, we shall call $(\mu_t)_{t\ge0}$ a {\it
differentiable} SC-semigroup if it is given by (\ref{3.2}). We
shall discuss how to extend a general SC-semigroup on $H$ to a
differentiable one on the entrance space $\tilde H$. For any
strongly continuous linear semigroup $(T_t)_{t\ge0}$ on $H$,
 \beqlb\label{3.3}
(\tilde T_t\tilde x)(s) = \tilde x(t+s),
\quad s,t>0,
 \eeqlb
defines a semigroup of linear operators $(\tilde T_t)_{t\ge 0}$
on $\tilde H$. Since
 \beqnn
\|\tilde T_t\tilde x\|_\sim^2
=
\int_0^\infty e^{-2bs}\|\tilde x(t+s)\|^2 ds
\le
\|T_t\|^2\int_0^\infty e^{-2bs}\|\tilde x(s)\|^2 ds,
 \eeqnn
we have $\|\tilde T_t\|_\sim \le \|T_t\|$. Let $(\tilde
U_\alpha)_{\alpha > b_0}$ denote the resolvent of $(\tilde
T_t)_{t\ge 0}$ and let $\tilde A$ denote its generator with domain
$D(\tilde A) = \tilde U_\alpha \tilde H \subset \tilde H$.

\blemma\label{l3.1} Let $J$ be defined as in Lemma \ref{l2.3}.
Then $JT_tx = \tilde T_t Jx$ for all $t\ge0$ and $x\in H$ and
$(\tilde T_t)_{t\ge 0}$ is a strongly continuous semigroup of
linear operators on $\tilde H$. \elemma

{\it Proof.} For $t\ge0$ and $x\in H$ we have
 \beqnn
JT_tx
=
\{T_sT_tx: s>0\}
=
\{T_tT_sx: s>0\}
=
\tilde T_tJx,
 \eeqnn
giving the first assertion. By the proof of Lemma \ref{l2.3},
$\|\tilde T_t\tilde x - \tilde
x\|_\sim = \| J\tilde x(t) - \tilde x\|_\sim \to 0$ as $t\to 0$,
that is, $(\tilde T_t)_{t\ge 0}$ is strongly continuous.
\qed

\blemma\label{l3.2} We have $\tilde U_\alpha \tilde x = \{U_\alpha
\tilde x(s): s>0\}$ and $\tilde A \tilde U_\alpha \tilde x =
\{AU_\alpha x(s): s>0\}$ for all $\tilde x\in \tilde H$. \elemma

{\it Proof.} The first assertion follows as we observe that
 \beqnn
\tilde U_\alpha\tilde x(s)
=
\int_0^\infty e^{-\alpha s}\tilde T_t\tilde x(s)dt
=
\int_0^\infty e^{-\alpha s} T_t\tilde x(s)dt
=
U_\alpha\tilde x(s),
 \eeqnn
and the second follows from the equality $\tilde A\tilde U_\alpha
\tilde x
= \alpha \tilde U_\alpha \tilde x - \tilde x$.
\qed

\btheorem\label{t3.1} All centered SC-semigroups associated with
$(T_t)_{t\ge0}$ satisfying (\ref{2.1}) are differentiable if and
only if all its locally square integrable entrance paths are
closable. \etheorem

{\it Proof.} Suppose that all entrance paths $\tilde x \in \tilde
H$ are closable and $(\mu_t)_{t\ge0}$ is an SC-semigroup given by
(\ref{2.8}). To each $\tilde x \in \tilde H$ there corresponds
some $\tilde x(0)\in H$ such that $\tilde x(s) = T_s\tilde x(0)$
for all $s>0$. This element is apparently determined by $\tilde x$
uniquely. Letting $\nu_0$ be the image of $\lambda_0$ under the
map $\tilde x \mapsto \tilde x(0)$ we get (\ref{3.2}). Conversely,
if $\tilde x= \{\tilde x(s): s>0\} \in \tilde H$ is not closable,
then
 \beqnn
\hat\mu_t(a)
=
\exp\bigg\{-\frac{1}{2}\int_0^t \lan \tilde x(s),a\ran ^2 ds\bigg\},
\quad t\ge0, a\in H^*,
 \eeqnn
defines a non-differentiable SC-semigroup.
\qed

\btheorem\label{t3.2}
All entrance paths for $(\tilde T_t)_{t\ge 0}$ are closable.
\etheorem

{\it Proof.} Suppose that $\bar x = \{\bar x(u): u>0\}$ is an
entrance path for $(\tilde T_t)_{t\ge0}$, where each $\bar x(u) =
\{\bar x(u,s): s>0\} \in \tilde H$ is an entrance path for
$(T_t)_{t\ge0}$. Then we get
 \beqlb\label{3.4}
\{\bar x(u,r+s): s>0\}
=
(\tilde T_r\bar x)(u)
=
\bar x(r+u)
=
\{\bar x(r+u,s): s>0\},
 \eeqlb
where the first equality follows from (\ref{3.3}) and the
second one holds since $\bar x$ is an entrance path for
$(\tilde T_t)_{t\ge0}$. Setting $\bar x(0) = \{\bar
x(s/2,s/2): s>0\}$ we have
 \beqlb\label{3.5}
\tilde T_u\bar x(0)(s)
=
\bar x(s/2,u+s/2)
=
\bar x(u,s),
 \eeqlb
where the first equality follows from (\ref{3.3}) and the
second one holds by (\ref{3.4}). Thus $\tilde T_u\bar x(0)
= \bar x(u)$, that is, $\bar x = \{\bar x(u): u>0\}$ is
closed by $\bar x(0)$. \qed

For the infinitely divisible probability measure $\lambda_0$ on
$\tilde H$ given by Theorem \ref{t2.2} we have
 \beqlb\label{3.6}
\hat\lambda_0(\tilde a)
=
e^{-\tilde \psi_0(\tilde a)},
\quad \tilde a\in \tilde H^*,
 \eeqlb
for a functional $\tilde \psi_0$ on $\tilde H^*$ with representation
 \beqlb\label{3.7}
\tilde \psi_0(\tilde a)
=
\frac{1}{2}\lan\tilde Ra,a\ran_\sim
- \int_{\tilde H^\circ} \left(e^{i\lan\tilde x,\tilde a\ran_\sim}
- 1 - i\lan\tilde x,\tilde a\ran_\sim \right) \tilde M(d\tilde x),
\quad \tilde a\in \tilde H^*,
 \eeqlb
where $R$ is a nuclear operator on $\tilde H$ and $\| \tilde
x\|_\sim^2 \tilde M(d\tilde x)$ is a finite measure on $\tilde
H^\circ := \tilde H \setminus \{0\}$; see e.g.\ \cite{L86}.

\btheorem\label{t3.3} Let $(\mu_t)_{t\ge 0}$ be a centered
SC-semigroup given by (\ref{2.8}) and let $\tilde\mu_t = J\mu_t$.
Then $(\tilde\mu_t)_{t\ge 0}$ is a differentiable centered
SC-semigroup associated with $(\tilde T_t)_{t\ge 0}$ and
 \beqlb\label{3.8}
\int_{\tilde H} e^{i\lan\tilde x,\tilde a\ran_\sim}
\tilde\mu_t(d\tilde x)
=
\exp \bigg\{-\int_0^t \tilde\psi_0(\tilde T_s^*\tilde a) ds\bigg\},
\quad t\ge0, \tilde a \in \tilde H^*.
 \eeqlb
\etheorem

{\it Proof.} By Lemma \ref{l2.3}, $J: H\mapsto \tilde H$ is an
embedding. Thus $(\tilde\mu_t)_{t\ge 0}$ is an SC-semigroup
associated with $(\tilde T_t)_{t\ge 0}$. Since $(T_t)_{t\ge 0}$ is
a strongly continuous semigroup, for any $\tilde a = \{\tilde
a(s): s>0\} \in \tilde H$ we have by dominated convergence
 \beqnn
& &\int_H \exp\bigg\{i\int_0^\infty e^{-2bs}\lan T_sx,\tilde a(s)
\ran ds\bigg\} \mu_t(d x)  \\
&=&
\lim_{n\to\infty}\int_H \exp\bigg\{i\sum_{k=1}^\infty n^{-1}
e^{-2bk/n}\lan T_{k/n}x,\tilde a(k/n)\ran \bigg\} \mu_t(d x)  \\
&=&
\lim_{n\to\infty}\int_H \exp\bigg\{i\bigg\lan x,\sum_{k=1}^\infty
n^{-1}e^{-2bk/n}T_{k/n}^*\tilde a(k/n)\bigg\ran \bigg\}\mu_t(d x)\\
&=&
\lim_{n\to\infty}\exp\bigg\{\int_0^t \bigg[\log\int_H \exp
\bigg\{i\bigg\lan x,\sum_{k=1}^\infty n^{-1}e^{-2bk/n}T_{k/n}^*
\tilde a(k/n)\bigg\ran \bigg\} \nu_s(d x) \bigg]ds\bigg\} \\
&=&
\lim_{n\to\infty}\exp\bigg\{\int_0^t \bigg[\log\int_H
\exp\bigg\{i\sum_{k=1}^\infty n^{-1}e^{-2bk/n} \lan T_{k/n}x,
\tilde a(k/n)\ran \bigg\} \nu_s(d x) \bigg]ds\bigg\}  \\
&=&
\exp\bigg\{\int_0^t \bigg[\log\int_H \exp\bigg\{i\int_0^\infty
e^{-2bu}\lan T_ux,\tilde a(u)\ran du\bigg\}
\nu_s(dx)\bigg] ds\bigg\}.
 \eeqnn
It follows that
 \beqlb\label{3.9}
\int_{\tilde H} e^{i\lan\tilde x,\tilde a\ran_\sim}
J\mu_t(d\tilde x)
=
\exp \bigg\{\int_0^t \bigg[\log\int_{\tilde H} e^{i\lan\tilde x,
\tilde a\ran_\sim} J\nu_s(d\tilde x)\bigg] ds\bigg\},
\quad t\ge0, \tilde a \in \tilde H.
 \eeqlb
Recall from the proof of Theorem \ref{t2.2} that $\nu_s$ is the
image of $\lambda_0$ under $\tilde x \mapsto \tilde x(s)$. Then
$\tilde T_s\lambda_0 = J\nu_s$ and (\ref{3.8}) follows from
(\ref{3.9}) and (\ref{3.6}). \qed

\btheorem\label{t3.4} Let $(\tilde \mu_t)_{t\ge0}$ be a centered
SC-semigroup associated with $(\tilde T_t)_{t\ge0}$ satisfying
 \beqlb\label{3.10}
\int_{\tilde H^\circ} \|\tilde x\|^2 \tilde\mu_t(d\tilde x)
< \infty,
\quad t\ge 0.
 \eeqlb
Then there is a centered SC-semigroup $(\mu_t)_{t\ge 0}$
associated with $(T_t)_{t\ge0}$ satisfying (\ref{2.1}) and
$\tilde\mu_t = J\mu_t$ for each $t\ge 0$. \etheorem

{\it Proof.} By Theorems \ref{t3.1} and \ref{t3.2}, $(\tilde
\mu_t)_{t\ge0}$ is differentiable, so it has the expression
(\ref{3.8}) for an infinitely divisible probability $\lambda_0$ on
$\tilde H$ defined by (\ref{3.6}). Then we get $(\mu_t)_{t\ge0}$
by Theorem \ref{t2.2}, which clearly satisfies the requirements.
\qed

By Theorems \ref{t3.3} and \ref{t3.4}, centered SC-semigroups
associated with $(T_t)_{t\ge0}$ and those associated with $(\tilde
T_t)_{t\ge0}$ are in 1-1 correspondence. Therefore we may reduce
some analysis of non-differentiable centered SC-semigroups to
those of differentiable ones studied in \cite{BRS96, FR00}.

\section{Ornstein-Uhlenbeck processes}
\setcounter{equation}{0}

In this section, we discuss constructions of the OU-processes. By
the results of the last section, a general centered SC-semigroup
on $H$ can be extended to a differentiable one on the entrance
space $\tilde H$. Then by Fuhrman and R\"ockner
\cite[Theorem~5.3]{FR00}, there is an extension $E$ of $\tilde H$
on which a c\`adl\`ag realization of the corresponding OU-process
can be constructed. In the sequel, we shall give a modification of
the arguments of Fuhrman and R\"ockner \cite{FR00} which provides
a smaller extension but still contains a c\`adl\`ag realization of
the OU-process. Fix $\alpha> b_0$ and define an inner product on
$\tilde H$ by
 \beqlb\label{4.1}
\lan\tilde x,\tilde y\ran_-
=
\int_0^\infty e^{-2bs}\lan U_\alpha\tilde x(s),U_\alpha
\tilde y(s)\ran ds,
\quad \tilde x,\tilde y\in \tilde H.
 \eeqlb
Let $\|\cdot\|_-$ be the corresponding norm and let $\bar H$ be the
completion of $\tilde H$ with respect to $\|\cdot\|_-$.

\blemma\label{l4.1} For $\tilde x\in \tilde H$ we have $\|\tilde
x\|_- \le \|U_\alpha\| \|\tilde x\|_\sim$, so the identity mapping
$I$ from $(\tilde H, \|\cdot\|_\sim)$ to $(\bar H,\|\cdot\|_-)$ is
a continuous embedding. \elemma

{\it Proof.} By (\ref{2.4}) and (\ref{4.1}),
 \beqnn
\|\tilde x\|_-^2
=
\int_0^\infty e^{-2bs} \|U_\alpha \tilde x(s)\|^2 ds
\le
\|U_\alpha\|^2\int_0^\infty e^{-2bs} \|\tilde x(s)\|^2 ds
=
\|U_\alpha\|^2\|\tilde x\|_\sim^2,
 \eeqnn
giving the desired estimate. \qed

Note that the embedding of $(\tilde H, \|\cdot\|_\sim)$ into
$(\bar H,\|\cdot\|_-)$ is not necessarily Hilbert-Schmidt, so our
extension is different from the one used in \cite{FR00}. For
$\tilde x\in \tilde H$ we have, by (\ref{4.1}),
 \beqnn
\|\tilde T_t\tilde x\|_-^2
=
\int_0^\infty e^{-2bs}\|U_\alpha T_t\tilde x(s)\|^2 ds
=
\int_0^\infty e^{-2bs}\|U_\alpha \tilde x(t+s)\|^2 ds
\le
e^{2bt}\|\tilde x\|_-^2.
 \eeqnn
Then each $\tilde T_t$ has a unique extension to a bounded linear
operator $\bar T_t$ on $\bar H$. Since the semigroup property and
strong continuity of $(\bar T_t)_{t\ge 0}$ hold on the dense
subspace $\tilde H$ of $\bar H$, they also hold on $\bar H$, that
is, the semigroup of linear operators $(\bar T_t)_{t\ge 0}$
extends $(\tilde T_t)_{t\ge 0}$. Let $(\bar
U_\alpha)_{\alpha>b_0}$ denote the resolvent of $(\bar T_t)_{t\ge
0}$ and let $\bar A$ denote its generator with domain $D(\bar A) =
\bar U_\alpha \bar H \subset \bar H$. Then $D(\bar A)$ is a
Hilbert space with inner product norm $\|\cdot\|_{\bar A}$ defined
by
 \beqlb\label{4.2}
\| \bar x\|_{\bar A}
=
\| \bar x\|_- + \| \bar A\bar x\|_-,
\quad \bar x\in D(\bar A).
 \eeqlb

\blemma\label{l4.2} We have $\tilde H\subset D(\bar A)$ and
 \beqlb\label{4.3}
\|\bar A\tilde y\|_-
\le
2(\alpha^2 \|U_\alpha\|^2 + 1)^{1/2} \|\tilde y\|_\sim,
\quad \tilde y\in \tilde H.
 \eeqlb
Consequently, the identity mapping $I$ from $(\tilde H,
\|\cdot\|_\sim)$ to $(D(\bar A),\|\cdot\|_{\bar A})$ is a
continuous embedding. \elemma

{\it Proof.} Suppose that $\tilde y \in D(\tilde A)\subset D(\bar
A)$. Then $\tilde y = \tilde U_\alpha \tilde x$ for some $\tilde
x\in \tilde H$. By (\ref{4.1}) and Lemma~\ref{l3.2},
 \beqnn
\|\bar A\tilde y\|_-^2
&=&
\| \tilde A \tilde U_\alpha \tilde x\|_-^2   \\
&=&
\int_0^\infty e^{-2bs} \|AU_\alpha^2 \tilde x(s)\|^2 ds  \\
&=&
\int_0^\infty e^{-2bs} \|\alpha U_\alpha^2 \tilde x(s)
- U_\alpha \tilde x(s)\|^2 ds  \\
&\le&
2(\alpha^2\|U_\alpha\|^2 + 1) \int_0^\infty e^{-2bs}
\|U_\alpha \tilde x(s)\|^2 ds  \\
&=&
2(\alpha^2\|U_\alpha\|^2 + 1) \|\tilde y\|_\sim^2.
 \eeqnn
Since $D(\tilde A)$ is a dense subset of $(\tilde H,
\|\cdot\|_\sim)$, we have $\tilde H\subset D(\bar A)$ by
(\ref{4.2}) and the above inequality remains true for all $\tilde
y\in \tilde H$. \qed

Now suppose that $(\mu_t)_{t\ge 0}$ and $(\tilde \mu_t)_{t\ge 0}$
are the SC-semigroups described in Theorem \ref{t3.3}. Let $\bar
\mu_t$ be the unique probability measure on $\bar H$ whose
restriction to $\tilde H$ is $\tilde \mu_t$. Then $(\bar
\mu_t)_{t\ge 0}$ is an SC-semigroup associated with $(\bar
T_t)_{t\ge 0}$. By (\ref{3.8}),
 \beqlb\label{4.4}
\int_{\bar H} e^{i\lan\bar x,\bar a\ran_-} \bar \mu_t(d\bar x)
=
\exp \bigg\{-\int_0^t \tilde\psi_0(\bar T_s^*\bar a) ds\bigg\},
\quad t\ge0, \bar a \in \bar H^* \subset \tilde H^*,
 \eeqlb
where $\tilde \psi_0(\cdot)$ is defined by (\ref{3.6}). Let
$(Q^{\bar\mu}_t)_{t\ge0}$ be the generalized Mehler semigroup
defined by (\ref{1.5}) from $(\bar T_t)_{t\ge0}$ and $(\bar
\mu_t)_{t\ge0}$. By \cite[Theorem~5.1]{FR00}, there is a
c\`adl\`ag $\tilde H$-valued process $\{\tilde Y_t: t\ge0\}$ with
$\tilde Y_0 =0$ and with independent increments such that $\tilde
Y_t - \tilde Y_r$ has distribution $\gamma_{t-r}$ with
 \beqlb\label{4.5}
\hat \gamma_{t-r}(\tilde a)
=
\exp\{-(t-r)\tilde\psi_0(\tilde a)\},
\quad t\ge r\ge 0, \tilde a\in \tilde H^*.
 \eeqlb
By the strong continuity of $(\bar T_t)_{t\ge0}$ and Lemma
\ref{l4.2}, $s\mapsto \bar T_{t-s} \bar A\tilde Y_s$ is a right
continuous $\bar H$-valued function of $s\in [0,t]$. Then for any
given $\bar x\in \bar H$ we may define the c\`adl\`ag $\bar
H$-valued process $\{\bar X_t: t\ge0\}$ by
 \beqlb\label{4.6}
\bar X_t
=
\bar T_t\bar x + \tilde Y_t + \int_0^t \bar T_{t-s}
\bar A\tilde Y_s ds,
\quad t\ge0.
 \eeqlb

\blemma\label{l4.3} The $\bar H$-valued random variable $\bar X_t$
has distribution $Q^{\bar\mu}_t(\bar x,\cdot)$ for every $t\ge0$.
In particular, if $\bar x\in JH$, then $\bar X_t\in JH$ a.s.\ for
every $t\ge0$. \elemma

{\it Proof.} We first prove that $\bar X_t^{(0)} := \bar X_t -
\bar T_t\bar x$ has distribution $\bar \mu_t(\cdot) = Q^{\bar
\mu}_t (0,\cdot)$. By the right continuity of $s\mapsto \bar
T_{t-s} \bar A \tilde Y_s$, we have
 \beqnn
\bar X^{(n)}_t
:= \tilde Y_t + \frac{t}{n}\sum_{k=1}^n \bar T_{(1-k/n)t}\bar A
\tilde Y_{kt/n} \to
\tilde Y_t + \int_0^t \bar T_{t-s}
\bar A\tilde Y_s ds
= \bar X_t^{(0)}
 \eeqnn
as $n\to \infty$. Let $D_0=0$ and $D_k = \bar T_{(1-n/n)t}\bar A
+ \cdots + \bar T_{(1-k/n)t}\bar A$. Then we have
 \beqnn
\bar X^{(n)}_t
&=&
\tilde Y_t + n^{-1}t[(D_1-D_2) \tilde Y_{t/n}
+ \cdots + (D_{n-1}-D_n) \tilde Y_{(n-1)t/n}
+ D_n\tilde Y_{nt/n}]  \\
&=&
(\tilde Y_{nt/n} - \tilde Y_{(n-1)t/n}) + \cdots
+ (\tilde Y_{2t/n} - \tilde Y_{t/n})
+ (\tilde Y_{t/n} - \tilde Y_0)     \\
& & + n^{-1}t[D_1(\tilde Y_{t/n} - \tilde Y_0)
+ D_2(\tilde Y_{2t/n} - \tilde Y_{t/n}) + \cdots   \\
& & + D_n(\tilde Y_{nt/n} - \tilde Y_{(n-1)t/n})]  \\
&=&
(I + n^{-1}tD_1)(\tilde Y_{t/n} - \tilde Y_0)
+ (I + n^{-1}tD_2)(\tilde Y_{2t/n} - \tilde Y_{t/n}) + \cdots   \\
& & + (I + n^{-1}tD_n)(\tilde Y_{nt/n} - \tilde Y_{(n-1)t/n})].
 \eeqnn
It follows that
 \beqnn
\E\exp\left\{i\lan\bar X^{(n)}_t,\bar a\ran_-\right\}
&=&
\E\exp\bigg\{i\sum_{k=1}^n\lan (I + n^{-1}tD_k)
(\tilde Y_{kt/n} - \tilde Y_{(k-1)t/n}),\bar a\ran_-\bigg\}  \\
&=&
\E\exp\bigg\{i\sum_{k=1}^n\lan (\tilde Y_{kt/n}
- \tilde Y_{(k-1)t/n}),
(I + n^{-1}tD_k)^*\bar a\ran_-\bigg\}  \\
&=&
\exp\bigg\{-\frac{t}{n}\sum_{k=1}^n
\tilde\psi_0((I + n^{-1}tD_k)^*\bar a)\bigg\}.
 \eeqnn
In view of (\ref{3.7}), $\tilde\psi_0(\cdot)$ is uniformly
continuous on any bounded subset of $\bar H^*$. Observe
also that
 \beqnn
& &\|\bar T^*_{(1-k/n)t}\bar a -
(I + n^{-1}tD_k)^*\bar a\|_-    \\
&=&
\bigg|\bigg| \bar T^*_{(1-k/n)t}\bar a - \bar a
- \frac{t}{n}\sum_{j=k}^n\bar T^*_{(1-j/n)t}
\bar A^*\bar a \bigg|\bigg|_{\bar H}  \\
&\le&
\sum_{j=k}^n\int_{(1-j/n)t}^{(1-(j-1)/n)t}\|\bar T^*_s
\bar A^*\bar a - \bar T^*_{(1-j/n)t} \bar A^*\bar a \|_- ds \\
&\le&
t \cdot \sup\left\{\|\bar T^*_{t_2}\bar A^*\bar a
- \bar T^*_{t_1} \bar A^*\bar a \|_-:
0\le t_1, t_2\le t \mbox{ and } |t_2-t_1| < t/n\right\},
 \eeqnn
which goes to zero as $n\to \infty$. Thus we have
 \beqnn
\E\exp\left\{i\lan\bar X_t^{(0)},\bar a\ran_-\right\}
&=&
\lim_{n\to \infty}\E\exp\left\{i\lan\bar X^{(n)}_t,
\bar a\ran_-\right\}    \\
&=&
\lim_{n\to \infty}\exp\bigg\{-\frac{t}{n}\sum_{k=1}^n
\tilde\psi_0(\bar T^*_{(1-k/n)t}\bar a)\bigg\}  \\
&=&
\exp\bigg\{-\int_0^t
\tilde\psi_0(\bar T^*_{t-s}\bar a)ds\bigg\},
 \eeqnn
so that $\bar X_t^{(0)}$ has distribution $Q^{\bar\mu}_t(0,\cdot)$.
Therefore, $\bar X_t$ has distribution $Q^{\bar\mu}_t(\bar x,\cdot)$.
If $\bar x= Jx$ for some $x\in H$, then $\bar T_t\bar x = \tilde
T_tJ x = J T_tx\in JH$ by Lemma {\ref{l3.1}. Since $\bar
\mu_t(\cdot)$ is supported by $JH$, so is $Q^{\bar\mu}_t(\bar
x,\cdot)$ and hence a.s.\ $\bar X_t\in JH$. \qed

\btheorem\label{t4.1} The process $\{\bar X_t: t\ge0\}$ defined by
(\ref{4.6}) is a c\`adl\`ag strong Markov process with transition
semigroup $(Q^{\bar\mu}_t)_{t\ge0}$. \etheorem

{\it Proof.} By the construction (\ref{4.6}), $\{\bar X_t:
t\ge0\}$ is adapted to the filtration ${\cal F}_t :=
\sigma(\{\tilde Y_s: 0\le s\le t\})$. For $r,t\ge0$,
 \beqnn
\bar X_{r+t} - \bar T_t\bar X_r
&=&
\tilde Y_{r+t} - \bar T_t\tilde Y_r
+ \int_r^{r+t} \bar T_{r+t-s} \bar A\tilde Y_s ds   \\
&=&
(\tilde Y_{r+t} - \tilde Y_r)
+ \int_r^{r+t} \bar T_{r+t-s} \bar A(\tilde Y_s - \tilde Y_r) ds.
 \eeqnn
Since $\{\tilde Y_{r+t} - \tilde Y_r: t\ge0\}$ given ${\cal F}_r$
is a process with independent increments and has the same law as
$\{\tilde Y_t: t\ge0\}$, Lemma \ref{l4.3} implies that
 \beqnn
\E\bigg[\exp\left\{i\lan\bar X_{r+t},\bar a\ran_-\right\}
\bigg|{\cal F}_r\bigg]
&=&
\exp\bigg\{i\lan\bar X_r,\bar T_t^*\bar a\ran_-
- \int_0^t \tilde\psi_0(\bar T_s^*\bar a) ds\bigg\}.
 \eeqnn
Thus $\{\bar X_t: t\ge0\}$ is a Markov process with transition
semigroup $(Q^{\bar\mu}_t)_{t\ge0}$. The strong Markov property
holds since $(Q^{\bar\mu}_t)_{t\ge0}$ is Feller. \qed

Now let $\bar x= Jx$ for some $x\in H$. In this case, $\bar X_t\in
JH$ a.s.\ by Lemma \ref{l4.3}. Recall that $J: H \to JH\subset
\tilde H\subset \bar H$ and let $X_t = 1_{JH}(\bar X_t) J^{-1}
(\bar X_t)$, where $J^{-1}: JH \to H$ denotes the inverse map of
$J$. Then $\{X_t: t\ge0\}$ is an OU-process with transition
semigroup $(Q^\mu_t)_{t\ge0}$ and $X_0=x$. This gives a
construction of the OU-process $\{X_t: t\ge0\}$ from the
c\`adl\`ag strong Markov process $\{\bar X_t: t\ge0\}$. In
general, $\{X_t: t\ge0\}$ does not have right continuous
modification in $H$. A similar construction in the measure-valued
setting has been used in \cite{L96} to prove the non-existence of
right continuous realization of a general immigration process.

\section{Brownian transition semigroups}
\setcounter{equation}{0}

We have seen that a general SC-semigroup on $H$ can always be
extended to a differentiable one in the entrance space $\tilde H$
and a c\`adl\`ag realization of the corresponding OU-process can
always be constructed in an extension $\bar H$ of $\tilde H$. In
this section, we provided some explicit characterization for the
non-negative elements of $\tilde L^2(\IR^d)$ and $\tilde
L^2(0,\infty)$ constructed respectively from $L^2(\IR^d)$ and
$L^2(0,\infty)$. It seems that the explicit characterization for
all elements of $\bar L^2(\IR^d)$ and $\bar L^2(0,\infty)$ is much
more sophisticated.

We first consider the case where $(T_t)_{t\ge0}$ is the transition
semigroup of the standard Brownian motion on $\IR^d$. Let $\tilde
H := \tilde L^2(\IR^d)$ be defined from $H:= L^2(\IR^d)$ and
$(T_t)_{t\ge0}$. Let
 \beqlb\label{5.1}
g_d(s,x) = \frac{1}{(2\pi s)^{d/2}}\exp\{-|x|^2/2s\},
\quad s>0, x\in \IR^d,
 \eeqlb
where $|\cdot|$ denote the Euclidean norm on $\IR^d$, and let
$S(\IR^d)$ be the set of signed-measures $\mu$ on $\IR^d$ with total
variation measures $\|\mu\|$ satisfying
 \beqlb\label{5.2}
\int_0^lds \int_{\IR^d}\|\mu\|(dx)\int_{\IR^d} g_d(2s,y-x)
\|\mu\|(dy)
< \infty
 \eeqlb
for some $l>0$. Let $S_+(\IR^d)$ and $\tilde L^2_+(\IR^d)$ denote
respectively the subsets of non-negative elements of $S(\IR^d)$
and $\tilde L^2(\IR^d)$.

\btheorem\label{t5.1} There is a 1-1 correspondence between
$\tilde x \in \tilde L^2_+(\IR^d)$ and $\mu \in S_+(\IR^d)$ which
is given by
 \beqlb\label{5.3}
\tilde x(s,\cdot) = \int_{\IR^d} g_d(s,\cdot -z) \mu(dz),
\quad s>0.
 \eeqlb
\etheorem

{\it Proof.} If $\mu \in S_+(\IR^d)$, then (\ref{5.3}) defines a
non-negative entrance path $\tilde x$ for $(T_t)_{t\ge0}$. Since
 \beqnn
\int_0^l \|\tilde x(s,\cdot)\|^2 ds
&=&
\int_0^lds \int_{\IR^d} \bigg(\int_{\IR^d} g_d(s,y-z)
\mu(dz)\bigg)^2 dy    \\
&=&
\int_0^lds \int_{\IR^d} dy \int_{\IR^d} g_d(s,y-x)\mu(dx)
\int_{\IR^d} g_s(y-z)\mu(dz) \\
&=&
\int_0^lds \int_{\IR^d}\mu(dx)\int_{\IR^d} g_d(2s,z-x) \mu(dz) \\
&<& \infty,
 \eeqnn
we have $\tilde x\in \tilde L^2_+(\IR^d)$. Conversely, suppose
that $\tilde x \in \tilde L^2_+(\IR^d)$ and let $\kappa_s(dy) =
\tilde x(s,y)dy$. Then $(\kappa_s)_{s>0}$ is a measure-valued
entrance law for $(T_t)_{t\ge0}$. By the property of the Brownian
semigroup, there is a measure $\mu$ on $\IR^d$ such that $\kappa_s
= \mu T_s$; see e.g.\ Dynkin \cite[p.80]{D89}. Thus $\tilde x(s,\cdot)$
has the representation (\ref{5.3}), and (\ref{5.2}) follows from
(\ref{2.2}) and the calculations above. \qed

When $d=1$, we can give a necessary and sufficient condition for
(\ref{5.2}). Observe that for $0<l\le 1$ we have
 \beqnn
\int_0^l g_1(2s,y-x) ds
<
\int_0^1 \frac{1}{2\sqrt{\pi s}} \exp\{-(y-x)^2/4\}ds
=
\frac{1}{\sqrt{\pi}} \exp\{-(y-x)^2/4\},
 \eeqnn
and for $l> 1$ we have
 \beqnn
\int_0^l g_1(2s,y-x) ds
>
\int_1^l \frac{1}{2\sqrt{\pi l}} \exp\{-(y-x)^2/4\}ds
=
\frac{l-1}{2\sqrt{\pi l}} \exp\{-(y-x)^2/4\}.
 \eeqnn
By Lemma \ref{l2.1} and the proof of Theorem \ref{t5.1},
(\ref{5.2}) holds if and only if
 \beqlb\label{5.4}
\int_{\IR}\|\mu\|(dx) \int_{\IR}\exp\{-(y-x)^2/4\} \|\mu\|(dy)
< \infty.
 \eeqlb

Theorem \ref{t5.1} gives a complete characterization of
non-negative elements of $\tilde L^2(\IR^d)$. By this result,
(\ref{5.3}) also defines an element of $\tilde L^2(\IR^d)$ for
$\mu\in S(\IR^d)$. Unfortunately, this representation does not
give all elements of $\tilde L^2(\IR^d)$. To see this, take any
sequence $\{a_k\} \subset \IR$ and observe that
 \beqnn
\int_0^\infty e^{-2bs} ds \int_{\IR} [g_1(s,y-x) - g_1(s,y)]^2 dy
\to 0
 \eeqnn
as $x\to0$. Then for each $k\ge1$ there exists $\varepsilon_k
\in (0,k^{-2})$ such that
 \beqlb\label{5.5}
a_k^2\int_0^\infty e^{-2bs} ds \int_{\IR} [g_1(s,y-\varepsilon_k)
- g_1(s,y)]^2 dy
\le
2^{-k}.
 \eeqlb
Let $x_k = k^{-1}$ and $z_k = k^{-1} + \varepsilon_k$. Then $z_k>
x_k> z_{k+1}> x_{k+1}> \dots$ decrease to zero. By (\ref{5.5}) and
the shift invariance of the Lebesgue measure it is easy to see
that
 \beqnn
\tilde x_n(s,\cdot)
=
\sum_{k=1}^n a_k[g_1(s,\cdot-z_k) - g_1(s,\cdot-x_k)],
\quad s>0,
 \eeqnn
defines a Cauchy sequence $\{\tilde x_n\} \subset \tilde L^2(\IR)$
with limit $\tilde x \in \tilde L^2(\IR)$ given by
 \beqlb\label{5.6}
\tilde x(s,\cdot)
=
\sum_{k=1}^\infty a_k[g_1(s,\cdot-z_k) - g_1(s,\cdot-x_k)],
\quad s>0.
 \eeqlb
To represent this element in the form of (\ref{5.3}) we need to
let
 \beqnn
\mu
=
\sum_{k=1}^\infty a_k\delta_{z_k}
-
\sum_{k=1}^\infty a_k\delta_{x_k},
 \eeqnn
which is clearly not belonging to $S(\IR)$ in general.

Now we consider the case where $(T_t)_{t\ge0}$ is the transition
semigroup of the absorbing barrier Brownian motion in $D=
(0,\infty)$. Let $\gamma(dy) = (1-e^{-y^2}) dy$ and let $\tilde
L^2(D,\gamma)$ be defined from $L^2(D,\gamma)$ and
$(T_t)_{t\ge0}$. Let
 \beqlb\label{5.7}
p_s(x,y) = g_1(s,y-x) - g_1(s,y+x),
\quad s,x,y>0,
 \eeqlb
and let
 \beqlb\label{5.8}
k_s(y) = 2^{-1}(d/dx) p_s(x,y)|_{x=0^+} = yg_1(s,y)/s,
\quad s,y>0.
 \eeqlb
Let $S(D,\gamma)$ be the set of signed-measures $\mu$ on $D$ with
total variation measures $|\mu|$ satisfying
 \beqlb\label{5.9}
\int_0^lds \int_D \bigg(\int_D p_s(x,y)|\mu|(dx)\bigg)^2
\gamma(dy) < \infty
 \eeqlb
for some $l>0$. Let $S_+(D,\gamma)$ and $\tilde L^2_+(D,\gamma)$
denote respectively the subsets of non-negative elements of
$S(D,\gamma)$ and $\tilde L^2(D,\gamma)$.

\btheorem\label{t5.2} There is a 1-1 correspondence between
$\tilde x \in \tilde L^2_+(D,\gamma)$ and $(a,\mu) \in [0,\infty)
\times S_+(D,\gamma)$ which is given by
 \beqlb\label{5.10}
\tilde x(s,\cdot) = ak_s(\cdot) + \int_D p_s(z,\cdot) \mu(dz),
\quad s>0.
 \eeqlb
\etheorem

{\it Proof.} If $(a,\mu) \in [0,\infty) \times S_+(D,\gamma)$, then
(\ref{5.10}) defines an entrance path $\tilde x$ for $(T_t)_{t\ge0}$.
Since
 \beqnn
\int_0^lk_s(y)^2 ds
\le
\int_0^\infty \frac{y^2}{2\pi s^3}e^{-y^2/s}ds
=
\frac{1}{2\pi y^2},
 \eeqnn
we have
 \beqnn
\int_0^l \|\tilde x(s,\cdot)\|^2 ds
&=&
\int_0^lds \int_D \bigg(ak_s(y)
+ \int_D p_s(z,y)\mu(dz)\bigg)^2 \gamma(dy)  \\
&\le&
2a^2\int_0^lds \int_D k_s(y)^2 \gamma(dy)   \\
& &
+ 2\int_0^lds \int_D \bigg(\int_D p_s(z,y)\mu(dz)\bigg)^2
\gamma(dy)    \\
&<& \infty,
 \eeqnn
that is, $\tilde x\in \tilde L^2_+(D,\gamma)$. Conversely, suppose
that $\tilde x \in \tilde L^2_+(\IR)$ and let $\kappa_s(dy) =
\tilde x(s,y)dy$. Then $(\kappa_s)_{s>0}$ is a measure-valued
entrance law for $(T_t)_{t\ge0}$. By the property of the absorbing
barrier Brownian motion, there is a constant $a\ge 0$ and a
measure $\mu$ on $D$ such that (\ref{5.10}) holds; see e.g.\
\cite[Lemma~1.1]{LS95}. Since
 \beqnn
\int_0^lds \int_D \bigg(\int_D p_s(z,y)\mu(dz)\bigg)^2 \gamma(dy)
\le
\int_0^l \|\tilde x(s,\cdot)\|^2 ds,
< \infty.
 \eeqnn
we have $\mu \in S(D,\gamma)$.
\qed

By the general results of section~4, an OU-process associated with
the absorbing barrier Brownian motion in $D= (0,\infty)$ always
has c\`adl\`ag realization in $\bar L^2(D,\gamma) \supset \tilde
L^2(D,\gamma)$ defined from $L^2(D,\gamma)$. It was observed in
\cite{DLSS03} that in a special case the process also has
c\`adl\`ag realization in $S(D,\gamma)$.

\bigskip

{\bf Acknowledgement} \ We thank Klaus Fleischmann and a referee
for their careful reading of the manuscript and helpful comments.


\noindent
\begin{thebibliography}{99}


\bibitem{BRS96}
Bogachev, V.I., R\"ockner, M. and Schmuland, B., {\it Generalized
Mehler semigroups and applications}, Probab. Th. Rel. Fields {\bf
105} (1996), 193-225.

\bibitem{DLSS04}
Dawson, D.A., Li, Z.H., Schmuland, B. and Sun, W., {\it
Generalized Mehler semigroups and catalytic branching processes
with immigration}, Potential Anal. \textbf{21}, 75-97.

\bibitem{D89}
Dynkin, E.B.: {\it Three classes of infinite dimensional diffusion
processes}, J. Funct. Anal. {\bf 86} (1989), 75-110.

\bibitem{FR00}
Fuhrman, M. and R\"ockner, M., {\it Generalized Mehler semigroups:
The non-Gaussian case}, Potential Anal. {\bf 12} (2000), 1-47.

\bibitem{L95/6}
Li, Z.H., {\it Convolution semigroups associated with
measure-valued branching processes}, Chinese Science Bulletin
(Chinese Edition) {\bf 40} (1995), 2018-2021 / (English Edition)
{\bf 41} (1996), 276-280.

\bibitem{L96}
Li, Z.H., {\it Immigration structures associated with
Dawson-Watanabe superprocesses}, Stochastic Process. Appl. {\bf
62} (1996), 73-86.

\bibitem{L98}
Li, Z.H., {\it Immigration processes associated with branching
particle systems}, Adv. Appl. Probab. {\bf 30} (1998), 657-675.

\bibitem{L02}
Li, Z.H., {\it Skew convolution semigroups and related immigration
processes}, Theory Probab. Appl. {\bf 46} (2002), 274-296.

\bibitem{LS95}
Li, Z.H. and Shiga T., {\it Measure-valued branching diffusions:
immigrations,  excursions and limit theorems}, J. Math. Kyoto
Univ. {\bf 35} (1995), 233-274.

\bibitem{L86}
Linde, W., {\it Probability in Banach Spaces -- Stable and
Infinitely Divisible Distributions}, Wiley, New York (1986).

\bibitem{N00}
van Neerven, J. M. A. M., {\it Continuity and representation of
Gaussian Mehler semigroups}, Potential Anal. {\bf 13} (2000),
199-211.

\bibitem{SS01}
Schmuland, B. and Sun W., {\it On the equation $\mu_{t+s} = \mu_s
* T_s\mu_t$}, Statist. Probab. Letters, {\bf 52} (2001), 183-188.

\bibitem{S88}
Sharpe, M.J., {\it General Theory of Markov Processes}, Academic
Press, New York (1988).



\end{thebibliography}
\end{document}